\documentclass[12pt]{amsart}

\usepackage{amssymb}
\usepackage{latexsym}
\usepackage{amsmath}
\usepackage{amscd}

\newtheorem{lemma}{Lemma}
\newtheorem{theorem}{Theorem}
\newtheorem{corollary}{Corollary}[section]

\usepackage{amscd, amsmath, amssymb, fancybox, eepic, graphicx}

\newcommand{\C}{\mathbb{C}}
\newcommand{\Q}{\mathbb{Q}}

\newcommand{\E}{\mathbb{E}}

\newcommand{\ch}{\operatorname{ch}}

\newcommand{\Tr}{\operatorname{Tr}}
\newcommand{\Wg}{\operatorname{Wg}}

\newcommand{\s}{\sigma}
\renewcommand{\t}{\tau}
\renewcommand{\a}{\alpha}
\renewcommand{\b}{\beta}
\renewcommand{\u}{\overline{u}}
\newcommand{\Sc}{\operatorname{Sc}}

\begin{document}

\title{Truncations of a Random Unitary Matrix and Young Tableaux}

\author{Jonathan Novak}
\address{Department of Mathematics and Statistics, Queen's University, Kingston,
Ontario, Canada}
\email{jnovak@mast.queensu.ca}
\thanks{}
\date{\today}
\begin{abstract}
Let $U$ be a matrix chosen randomly, with respect 
to Haar measure, from the unitary group $U(d).$  We express the 
moments of the trace of any submatrix of $U$ as a 
sum over partitions whose terms count 
certain standard and semistandard Young tableaux.
Using this combinatorial interpretation, we obtain a simple
closed form for the moments of an individual entry of a random
unitary matrix and use this to deduce that the rescaled entries 
converge in moments
to standard complex Gaussian random variables.  In addition, we recover a 
well-known theorem of E. Rains which shows that the moments of the trace
of a random unitary matrix enumerate permutations with restricted increasing
subsequence length.
\end{abstract}
\maketitle

\section{Introduction}
Consider the unitary group $U(d)$ as a probability space under normalized
Haar measure.  Given a random variable $X:U(d) \rightarrow \C,$ 
its expected value is defined to be
	$$\E_{U(d)}(X)=\int_{U(d)} X dU.$$
When studying a random variable $X,$ one often wishes to know its
moments
	$$\E_{U(d)}(X^m \overline{X}^n),$$
since in many situations the moments of $X$ uniquely determine its
distribution.
$X$ is called a polynomial random variable if 
there is a polynomial $f \in \C[x_{11}, \dots, x_{dd}]$ such that
	$$X(U) = f(u_{11}, \dots, u_{dd})$$ 
for all $U \in U(d).$  If $X$ is a polynomial random variable,
then
	$$\E_{U(d)}(X^m \overline{X}^n) = 0\ \text{if } m \neq n,$$
(see \cite{CS}), so knowledge of the moments of $X$ reduces to 
knowledge of the quantities
	$$E_{U(d)}(|X|^{2n}).$$

\smallskip

Recently, certain polynomial random variables on the unitary group
have been shown to possess interesting combinatorial properties. 
For a random matrix $U \in U(d),$ write its characteristic polynomial
as
	$$\det(U-zI) = (-1)^d \sum_{j=0}^d (-1)^j \Sc_j(U) z^{d-j}.$$
$\Sc_j(U)$ is called the $j$-th {\it secular coefficient} of $U.$
In particular, 
	\begin{align*}
	\Sc_1(U) &= \Tr(U) \\
	\Sc_d(U) &= \det(U).
	\end{align*}	
It was shown by Rains in \cite{Rains} that 
	$$\E_{U(d)}(|\Tr(U)|^{2n}) = \sum_{\lambda \vdash n, \ell(\lambda)
	\leq d} (f^{\lambda})^2,$$
the number of pairs of semistandard Young tableaux on the same Young
diagram, over all possible Young diagrams with $n$ boxes and at most 
$d$ rows.  By the Schensted correspondence, this is equal to the number
of permutations in the symmetric group $S_n$ with no increasing subsequence
of length greater than $d.$  

\smallskip

The other secular coefficients of a random unitary matrix also
encode interesting combinatorial information, as shown by Diaconis
and Gamburd in \cite{DG}.  For $jn \leq d,$ \cite{DG} proves that
	$$\E_{U(d)}(|\Sc_j(U)|^{2n}) = H_n(j),$$
where $H_n(j)$ is the number of $n \times n$ matrices whose entries
are nonnegative integers and whose rows and columns all sum to $j$ 
(``magic squares'').  The method of proof given
in \cite{DG} is simple and elegant, and consists of two main ingredients:
	\begin{itemize}
		\item
		The secular coefficients of $U$ are the elementary
		symmetric functions applied to the eigenvalues of 
		$U.$
		
		\item
		The elementary symmetric functions can be written
		as linear combinations of Schur functions, which are 
		the irreducible characters of the 
		unitary group and thus satisfy orthogonality relations. 
	\end{itemize}

\smallskip

Now suppose that we are given a random matrix $U \in U(d),$ and we 
want to calculate the moments of a single entry of $U$
	$$\E_{U(d)}(|u_{ij}|^{2n})$$
(note that since the permutation matrices are in $U(d)$ and Haar measure
is translation invariant, all entries of $U$ are equidistributed).
More generally, for a positive integer $k$ with $1 \leq k \leq d,$
let $U_k$ denoted the $k \times k$ upper left corner of $U.$
We could ask about the moments of the matrix $U_k$
	$$\E_{U(d)}(|\Tr(U_k)|^{2n})$$
(again, since the permutation matrices are in $U(d),$ the traces
of any two $k \times k$ submatrices of $U$ are equidistributed).
The methods used in \cite{Rains}, \cite{DG} will not work in this 
situation, since we are no longer dealing with the eigenvalues of
a unitary matrix, but rather the eigenvalues of submatrices of a 
random unitary matrix.
However, quite surprisingly, there is a simple
combinatorial formula for these moments.

\begin{theorem}
Let $U$ be a matrix chosen randomly with respect to Haar measure
from the unitary group $U(d),$  and let
$U_k$ be its $k \times k$ upper left corner.  
We have	
	$$\E_{U(d)}(|\Tr(U_k)|^{2n}) = \sum_{\lambda \vdash n,\
        \ell(\lambda) \leq k} (f^{\lambda})^2 
	\frac{s_{\lambda,k}(1)}{s_{\lambda,d}(1)}.$$
Here, for a positive integer $r,$ $s_{\lambda,r}(1)$ denotes the number
of semistandard Young tableaux on the shape $\lambda$ with entries
from the set $[r]=\{1, \dots, r\}.$
\end{theorem}
	
\smallskip

The notation $s_{\lambda,r}(1)$ is shorthand for the value of the 
Schur function $s_{\lambda}$ obtained by setting the first $r$ variables
equal to $1$ and making the remaining variables $0$ 
	$$s_{\lambda}(1,1, \dots, 1, 0, 0 ,\dots).$$

\smallskip

The matrix $U_k$ is called a {\it truncation} of $U.$  Truncations of
random unitary matrices where first studied in \cite{SZ} from an
analytical point of view, where it was shown that $U_k$ is a contraction
(i.e. all of its eigenvalues lie in the closed unit disc in $\C$), and
the joint probability density function of the eigenvalue was found to
be
	$$C_{d,k} \prod_{i<j} |z_i-z_j|^2 \prod_{j=1}^k (1-|z_j|^2)^{d-k-1}.$$
$C_{d,k}$ is a normalization constant that was found in \cite{PR} to 
be
	$$C_{d,k} = \frac{1}{\pi^k k!} \prod_{j=0}^{k-1}
	{d-k+j-1 \choose j} (d-k+j).$$
	
\smallskip

Note that Rains's result is an immediate corollary of Theorem 1; it
is simply the case $k=d.$

\begin{corollary}
In the special case $k=d,$ we have 
$$\E_{U(d)}(|\Tr(U)|^{2n}) = \sum_{\lambda \vdash n,
        \ell(\lambda) \leq d}
        (f^{\lambda})^2.$$
\end{corollary}

\begin{proof}
When $k=d,$
        $$\frac{s_{\lambda,k}(1)}{s_{\lambda,d}(1)}
        =\frac{s_{\lambda,d}(1)}{s_{\lambda,d}(1)}
        =1.$$
\end{proof}
\smallskip

Rains's theorem is the extreme case $k=d$ of Theorem 1.  The 
other extreme $k=1$ is also of interest, since this corresponds to the
computation of the moments of a single entry of a random unitary matrix.

\begin{corollary}
In the special case $k=1,$ we have
	$$\E_{U(d)}(|u_{ij}|^{2n}) 
	= \frac{n!} {d(d+1) \dots (d+n-1)} = {d+n-1 \choose n}^{-1}.$$
\end{corollary}

\begin{proof}
For $k=1$ the only contribution to the sum is made by the single
partition whose diagram is a row of $n$ boxes.
Thus we have
	$$\int_{U(d)} |u_{ij}|^{2n}dU
        = \frac{1} {s_{n,d}(1)}.$$
The generalized hook length formula asserts that
	$$s_{\lambda,d}(1) = \prod_{\square \in \lambda}
	\frac{d+c(\square)}{h(\square)},$$
where $c(\square)$ is the content of the box, and $h(\square)$ is its
hook length (see \cite{Stanley}).  For the single row partition of $n,$
this gives
	$$s_{n,d}(1) = \frac{d(d+1) \dots (d+n-1)}{n!},$$
and the result follows.
\end{proof}

\smallskip

Explicitly knowing the moments of $u_{ij}$ makes it easy to determine
its limiting distribution.
Recall that if $x,y$ are Gaussian random variables with mean $0$ and
variance $1/2,$ then the random variable $z=x+iy$ is called a standard
complex Gaussian.

\begin{corollary}
As $d \rightarrow \infty,$ the random variable $\sqrt{d}u_{ij}$ converges
in distribution to a standard complex Gaussian random variable.
\end{corollary}

\begin{proof}
It is well-known (see for instance \cite{PR}) that the 
moments of a standard complex Gaussian $z$ are given by 
	$$\E(z^m \overline{z}^n) = \delta_{mn} n!.$$
Corollary 1.2 shows that
	\begin{align*}
	\E_{U(d)}((\sqrt{d}u_{ij})^m (\sqrt{d}\u_{ij})^n) &=
	\delta_{mn} \frac{d^n n!}{d(d+1) \dots (d+n-1)} \\
	&= \delta_{mn} \frac{n!}{(1)(1+\frac{1}{d}) \dots (1+\frac{n-1}{d})} \\
	& \rightarrow \delta_{mn} n!.
	\end{align*}
\end{proof}

\smallskip

In order to prove Theorem 1, one needs to connect unitary expectations
to symmetric function theory by some method other than applying symmetric
functions to eigenvalues.  This can be done using the Weingarten function
introduced in \cite{CS}, which is a powerful tool for computing
the moments of polynomial random variables on the unitary group.  The
Weingarten function has already been used in free probability theory
to prove asymptotic freeness results for random unitary matrices
(see the recent book \cite{NS} for a clear account of this).  

\section{The Weingarten Function}

For any positive integers $d$ and $n,$ define a function
$\Wg(d,n,\cdot):S_n \rightarrow \Q$ by
$$\Wg(d,n,\s):= \frac{1}{n!^2} \sum_{\lambda \vdash n,\ \ell(\lambda) \leq d}
\frac{(f^{\lambda})^2}{s_{\lambda,d}(1)} \chi^{\lambda}(\s),$$
where $\chi^{\lambda}$ is the irreducible character
of $S_n$ labeled by $\lambda.$

The following integration formula was proved in \cite{CS}. 

\begin{theorem}
Let $i,j,i',j':[n] \rightarrow [d]$ be any functions.  Then
\begin{align*}
& \int_{U(d)}u_{i(1)j(1)} \dots u_{i(n)j(n)}
\u_{i(1)j(1)} \dots \u_{i(n)j(n)}dU \\ &=
\sum_{\s,\t \in S_n} \delta_{i(1)i'(\s(1))} \dots \delta_{i(n)i'(\s(n))}
\delta_{j(1)j'(\t(1))} \dots \delta_{j(n)j'(\t(n))}
\Wg(d,n,\t \s^{-1}),
\end{align*}
where $\delta$ is the Kronecker delta.
\end{theorem}

\smallskip

We can succinctly express the moments of $\Tr(U_k)$ in terms
of the Weingarten function as follows:

\begin{lemma}
For any positive integers $n$ and $k,$ where $1 \leq k \leq d,$ we have
$$\E_{U(d)}(|\Tr(U_k)|^{2n}) = n! 
\sum_{\a}{n \choose \a} \sum_{\s \in S_{\a}}
\Wg(d,n,\s),$$
where the outer sum runs over all weak $k$-part compositions $\a$ of
$n,$ and the inner sum runs over all permutations in the Young subgroup
$S_{\a}$ of $S_n.$
\end{lemma}

\begin{proof}
This is really just a calculation.  We expand

	\begin{align*}
	|\Tr(U_k)|^{2n} &= |u_{11}+ \dots + u_{kk}|^{2n} \\
	&= (u_{11}+ \dots + u_{kk})^n (\u_{11}+ \dots + \u_{kk})^n \\
	&= \sum_{\a} \sum_{\b} {n \choose \a}{n \choose \b}
	u^{\a}\u^{\b},
	\end{align*}

where we are summing over all pairs of weak $k$-part compositions of 
$n,$ $\a = (a_1, \dots, a_k)$ and $\b = (b_1, \dots, b_k).$
We are using multi-index notation,

	\begin{align*}
	{n \choose \a} &= {n \choose a_1, \dots, a_k} \\
	u^{\a} &= u_{11}^{a_1} \dots u_{kk}^{a_k} \\
	{n \choose \b} &= {n \choose b_1, \dots, b_k} \\
        \u^{\b} &= \u_{11}^{b_1} \dots \u_{kk}^{b_k} 
	\end{align*}
	
Hence 
		$$\E_{U(d)}(|\Tr(U_k)|^{2n}) = 
		\sum_{\a} \sum_{\b} {n \choose \a}{n \choose \b}
        	\E_{U(d)}(u^{\a}\u^{\b}).$$ 
We will use the Weingarten integration formula to evaluate the expectation
$\E_{U(d)}(u^{\a}\u^{\b})$ for a fixed pair of compositions $\a,\b.$
Implicitly define coordinate functions $i_{\a},j_{\a},i_{\b},j_{\b}
:[n] \rightarrow [d]$ by setting
\begin{align*}
	&\int_{U(d)} u_{i_{\a}(1)j_{\a}(1)} \dots 
	u_{i_{\a}(n)j_{\a}(n)} 
	\u_{i_{\b}(1)j_{\b}(1)} \dots \u_{i_{\b}(n)j_{\b}(n)} dU \\
	&:= \int_{U(d)}u_{11}^{a_1} \dots u_{kk}^{a_k}
	\u_{11}^{b_1} \dots \u_{kk}^{b_k} dU \\
	&= \E_{U(d)}(u^{\a}\u^{\b}).
\end{align*}
Applying the Weingarten integration formula, we have
\begin{align*}
	&\E_{U(d)}(u^{\a}\u^{\b}) \\
	&= \sum_{\s,\t \in S_n} 
	\delta_{i_{\a}(1)i_{\b}(\s(1))} \dots \delta_{i_{\a}(n)i_{\b}(\s(n))}
	\delta_{j_{\a}(1)j_{\b}(\t(1))} \dots \delta_{j_{\a}(n)j_{\b}(\t(n))}
	\Wg(d,n,\t \s^{-1}).
\end{align*} 
Since we are only taking entries from the diagonal, we have 
$i_{\a}=j_{\a}$ and $i_{\b}=j_{\b}.$  Moreover, the level sets of these
functions are easy to read off:
	\begin{align*}
	i_{\a}^{-1}(1) &= [1,a_1]\\
	i_{\a}^{-1}(2) &= [a_1+1,a_1+a_2] \\ 
	&\vdots		\\
	i_{\a}^{-1}(k) &= [n-a_k+1, n]
	\end{align*}
and
	\begin{align*}
        i_{\b}^{-1}(1) &= [1,b_1]\\
        i_{\b}^{-1}(2) &= [b_1+1,b_1+b_2] \\
        &\vdots         \\
        i_{\b}^{-1}(k) &= [n-b_k+1, n]
        \end{align*}
Hence
	\begin{align*}
	\prod_{r=1}^{a_1}\delta_{i_{\a}(r)i_{\b}(\s(r))}
	&= \prod_{r=1}^{a_1}\delta_{1i_{\b}(\s(r))} \\
	\prod_{r=a_1+1}^{a_1+a_2}\delta_{i_{\a}(r)i_{\b}(\s(r))}
        &= \prod_{r=a_1+1}^{a_1+a_2}\delta_{2i_{\b}(\s(r))} \\
	& \vdots \\
	\prod_{r=n-a_k+1}^n \delta_{i_{\a}(r)i_{\b}(\s(r))}
        &= \prod_{r=n-a_k+1}^n \delta_{ki_{\b}(\s(r))}, 
	\end{align*}
Thus in order for the product
$$\delta_{i_{\a}(1)i_{\b}(\s(1))} \dots \delta_{i_{\a}(n)i_{\b}(\s(n))}$$
to be nonzero, we see that $\s$ must bijectively map the interval
$[1,a_1]$ onto the interval $[1,b_1],$ and $\s$ must also bijectively map the 
interval $[a_1+1,a_1+a_2]$ onto the interval $[b_1+1,b_1+b_2],$ etc.
Similarly, in order for the product
$$\delta_{j_{\a}(1)j_{\b}(\t(1))} \dots \delta_{j_{\a}(n)j_{\b}(\t(n))}$$
to be nonzero, we see that $\t$ must bijectively map the interval
$[1,a_1]$ onto the interval $[1,b_1],$ and $\s$ must also bijectively map the
interval $[a_1+1,a_1+a_2]$ onto the interval $[b_1+1,b_1+b_2],$ etc.
Thus we see that the expectation $\E_{U(d)}(u^{\a}\u^{\b})$ is zero
unless:
	\begin{itemize}
		\item
		$\a = \b,$ i.e. these two are the same weak $k$-part
		composition of $n,$

		\item
		$\s,\t$ are both in the Young subgroup $S_{\a},$
		i.e. the subgroup of permutations in $S_n$ that
		permute the first $a_1$ symbols amongst themselves,
		the next $a_2$ symbols amongst themselves, etc.
	\end{itemize}
Thus 
	\begin{align*}
	\E_{U(d)}(|\Tr(U_k)|^{2n}) &= 
                \sum_{\a} {n \choose \a}^2 \E_{U(d)}(u^{\a}\u^{\a}) \\
		&= \sum_{\a} {n \choose \a}^2 \sum_{\s,\t \in S_{\a}}
		\Wg(d,n,\t \s^{-1}) \\
		&= \sum_{\a} {n \choose \a}^2 \alpha! \sum_{\s \in S_{\a}}
                \Wg(d,n,\s) \\ 		
		&= n! \sum_{\a} {n \choose \a} \sum_{\s \in S_{\a}}
                \Wg(d,n,\s) 
	\end{align*}
\end{proof}

\smallskip

\section{Proof of the Main Theorem}

We are now in a position to prove Theorem 1.  

\begin{proof}
Having proved that
	$$\E_{U(d)}(|\Tr(U_k)|^{2n}) = n! 
	\sum_{\a} {n \choose \a} \sum_{\s \in S_{\a}}
                \Wg(d,n,\s),$$
we will work with the sum on the right.  Plugging in the definition of the
Weingarten function, we have
	$$n! \sum_{\a} {n \choose \a} \sum_{\s \in S_{\a}}
                \Wg(d,n,\s) 
	= n! \sum_{\a} {n \choose \a} \sum_{\s \in S_{\a}}
	\frac{1}{n!^2} \sum_{\lambda \vdash n,\ \ell(\lambda) \leq d}
\frac{(f^{\lambda})^2}{s_{\lambda,d}(1)} \chi^{\lambda}(\s).$$	
Changing order of summation, this becomes
	\begin{align*}
	\sum_{\lambda \vdash n,\ \ell(\lambda) \leq d} (f^{\lambda})^2
	\frac{1}{s_{\lambda,d}(1)} \sum_{\a} \frac{1}{\a !}
	\sum_{\s \in S_{\a}} \chi^{\lambda}(\s) 
	&= \sum_{\lambda \vdash n,\ \ell(\lambda) \leq d} (f^{\lambda})^2
        \frac{1}{s_{\lambda,d}(1)} \sum_{\a} \langle 1, \chi^{\lambda} 
	\rangle_{S_{\a}},
	\end{align*}
where the inner product $\langle \cdot, \cdot \rangle_{S_{\a}}$ is the 
averaged dot product on the space $CF(S_{\a})$ of complex-valued class
functions on the group $S_{\a}.$
Note that this sum may be written as
$$\sum_{\lambda \vdash n,\ \ell(\lambda) \leq d} (f^{\lambda})^2
        \frac{1}{s_{\lambda,d}(1)}
	\sum_{\alpha} \langle 1,
                \chi^{\lambda}\downarrow^{S_n}_{S_{\alpha}}
                \rangle_{S_{\alpha}},$$
where $\chi^{\lambda}\downarrow^{S_n}_{S_{\alpha}}$ is the restriction
of the irreducible character $\chi^{\lambda}$ of $S_n$ to the subgroup
$S_{\a}.$  Now, the function which is identically $1$ is the character of 
the trivial representation of $S_{\a}.$  Thus we may apply Frobenius 
reciprocity (see for instance \cite{Sagan});
	$$\langle 1,
                \chi^{\lambda}\downarrow^{S_n}_{S_{\alpha}} \rangle_{S_{\a}}
	= \langle 1\uparrow_{S_{\alpha}}^{S_n},
                \chi^{\lambda}
                 \rangle_{S_n},$$
where $1\uparrow_{S_{\alpha}}^{S_n}$ is the induction of the 
trivial character of $S_{\a}$ to $S_n.$  

\smallskip

The final step in the proof relies on the characteristic map
$\ch^n:CF(S_n) \rightarrow \Lambda^n,$ where $\Lambda^n$ is the 
inner product space of degree $n$ symmetric functions equipped with the
Hall inner product $\langle \cdot, \cdot \rangle_{\Lambda^n}$ 
(see for \cite{Sagan} or \cite{Stanley}).  
The class function is an isometry, and has the following important
properties: 
	\begin{align*}
		\ch_n(1\uparrow_{S_{\alpha}}^{S_n}) &= h_{\a} \\
		\ch_n(\chi^{\lambda}) &= s_{\lambda},
	\end{align*}
where $h_{\a}$ is the complete homogeneous symmetric function
indexed by $\a,$ and $s_{\lambda}$ is the Schur function indexed
by $\lambda.$
Thus we have
	$$\langle 1\uparrow_{S_{\alpha}}^{S_n},
                \chi^{\lambda}
                 \rangle_{S_n} = 
		\langle \ch_n(1\uparrow_{S_{\alpha}}^{S_n}),
                \ch_n(\chi^{\lambda})
                 \rangle_{S_n}
		= \langle h_{\a}, s_{\lambda} \rangle_{\Lambda^n}.$$
It is well-known that the Schur functions constitute an
orthonormal basis for 
$\Lambda^n,$ and that the coordinates of the complete homogeneous
symmetric functions with respect to the basis of Schur functions
are the Kostka numbers (see \cite{Stanley}).  That is,
	$$h_{\a} = \sum_{\mu \vdash n} K_{\mu \a}s_{\mu},$$
where the Kostka number $K_{\mu \a}$ is by definition the number of 
semistandard Young tableaux on the diagram of $\mu$ with content vector
$\a.$  Thus,
	$$\langle h_{\a}, s_{\lambda} \rangle_{\Lambda^n}
	= \sum_{\mu \vdash n} K_{\mu \a} 
	\langle s_{\mu}, s_{\lambda} \rangle_{\Lambda^n}
	= K_{\lambda, \a}.$$

\smallskip

Thus we have
	\begin{align*}
		\sum_{\lambda \vdash n,\ \ell(\lambda) \leq d} (f^{\lambda})^2
        \frac{1}{s_{\lambda,d}(1)} \sum_{\a} \langle 1, \chi^{\lambda}
        \rangle_{S_{\a}} &=
	\sum_{\lambda \vdash n,\ \ell(\lambda) \leq d} (f^{\lambda})^2
        \frac{1}{s_{\lambda,d}(1)} \sum_{\a} K_{\lambda \a} \\
	&= \sum_{\lambda \vdash n,\ \ell(\lambda) \leq d} (f^{\lambda})^2
        \frac{s_{\lambda,k}(1)}{s_{\lambda,d}(1)},
	\end{align*}
where the last equality follows from the fact that, by definition,
	$$\sum_{\a}K_{\lambda \a} = s_{\lambda,k}(1),$$
since the sum runs over all weak $k$-part compositions of $n.$ 

Finally, we remark that if $\lambda$ is a partition of $n$ with 
$\ell(\lambda) >k,$ then $s_{\lambda,k}(1)=0.$  Hence,
$$\sum_{\lambda \vdash n,\ \ell(\lambda) \leq d} (f^{\lambda})^2
        \frac{s_{\lambda,k}(1)}{s_{\lambda,d}(1)} =
\sum_{\lambda \vdash n,\ \ell(\lambda) \leq k} (f^{\lambda})^2
        \frac{s_{\lambda,k}(1)}{s_{\lambda,d}(1)},$$
which proves our theorem.

\end{proof}

\section{Conclusion}
In this paper, we have only investigated the moments of the trace
of a truncation of a random unitary matrix.  It seems possible to analyze
the moments of the other secular coefficients of a truncation by the 
same method, and it would be interesting to see what combinatorial
interpretations can be given for the moments of these coefficients.

\smallskip

In \cite{DG}, the moments of secular coefficients of random orthogonal
and symplectic matrices are investigated, and results analogous to the 
unitary case are proved.  In \cite{CS}, a notion of Weingarten function
is defined for the orthogonal and symplectic groups.  It is likely 
possible to analyze the moments of secular coefficients of truncations of
random orthogonal and symplectic matrices using the Weingarten function
for these groups.

\section{Acknowledgements}
I am grateful to both Roland Speicher and Jamie Mingo for several
helpful discussions, and to Richard Stanley for some encouraging comments.

\end{document}